\documentclass{amsart}

\usepackage{amsmath}
\usepackage{amsfonts}
\usepackage{amssymb}

\numberwithin{equation}{section}

\newtheorem{theorem}[subsection]{Theorem}
\newtheorem{proposition}[subsection]{Proposition}
\newtheorem{corollary}[subsection]{Corollary}
\newtheorem{lemma}[subsection]{Lemma}
\newtheorem{claim}[subsubsection]{Claim}

\theoremstyle{definition}

\newcommand{\wh}{\widehat}

\title{An application of a local version of Chang's theorem}

\author{T. Sanders}

\begin{document}

\begin{abstract}
Suppose that $G$ is a compact Abelian group. If $A \subset G$ then
how small can $\|\chi_A\|_{A(G)}$ be? In general there is no
non-trivial lower bound.

In \cite{BJGSK} Green and Konyagin showed that if $\widehat{G}$
has sparse small subgroup structure and $A$ has density $\alpha$
with $\alpha(1-\alpha) \gg 1$ then $\|\chi_A\|_{A(G)}$ does admit
a non-trivial lower bound. To complement this \cite{TWSCVS}
addressed the case where $\wh{G}$ has rich small subgroup
structure and further claimed a result for general compact Abelian
groups. In this note we prove this claim by fusing the techniques
of \cite{BJGSK} and \cite{TWSCVS} in a straightforward fashion.
\end{abstract}

\maketitle

\section{Notation and introduction}\label{intro}

We use the Fourier transform on compact Abelian groups, the basics
of which may be found in Chapter 1 of Rudin \cite{WR}; we take a
moment to standardize our notation.

Suppose that $G$ is a compact Abelian group. Write $\widehat{G}$
for the dual group, that is the discrete Abelian group of
continuous homomorphisms $\gamma:G \rightarrow S^1$, where
$S^1:=\{z \in \mathbb{C}:|z|=1\}$. $G$ may be endowed with Haar
measure $\mu_G$ normalised so that $\mu_G(G)=1$ and as a
consequence we may define the Fourier transform
$\widehat{.}:L^1(G) \rightarrow \ell^\infty(\widehat{G})$ which
takes $f \in L^1(G)$ to
\begin{equation*}
\widehat{f}: \widehat{G} \rightarrow \mathbb{C}; \gamma \mapsto
\int_{x \in G}{f(x)\overline{\gamma(x)}d\mu_G(x)}.
\end{equation*}
We write
\begin{equation*}
A(G):=\{f \in L^1(G): \|\widehat{f}\|_1 < \infty\},
\end{equation*}
and define a norm on $A(G)$ by $\|f\|_{A(G)}:=\|\widehat{f}\|_1$.

The following result is well known.
\begin{proposition}\label{Cohen}
Suppose that $G$ is a compact Abelian group. Suppose that $A
\subset G$ has density $\alpha$ and for all finite $V \leq \wh{G}$
we have\footnote{Here $\{x\}$ denoted the fractional part of $x$.}
$\{\alpha|V|\}(1-\{\alpha|V|\})>0$. Then $\chi_A \not \in A(G)$.
\end{proposition}
In \cite{TWSCVS} the following quantitative version of the above
result was claimed.
\begin{theorem}\label{maintheorem}
Suppose that $G$ is a compact Abelian group. Suppose that $A
\subset G$ has density $\alpha$ and for all finite $V \leq \wh{G}$
with $|V| \leq M$ we have $\{\alpha|V|\}(1-\{\alpha|V|\})\gg 1$.
Then
\begin{equation*}
\|\chi_A\|_{A(G)} \gg \log \log \log M.
\end{equation*}
\end{theorem}
The objective of these notes is to prove Theorem
\ref{maintheorem}; for completeness and to illuminate the
arguments we begin with a proof of Proposition \ref{Cohen}. In
fact Proposition \ref{Cohen} follows is a straightforward
corollary of the celebrated idempotent theorem of Cohen \cite{PC}
and, indeed, the instance of the idempotent theorem which we use
was in fact proved even earlier by Rudin in \cite{WRID}.

\section{Proof of Proposition \ref{Cohen}}

We shall prove the following stronger result which characterizes
those sets $A \subset G$ for which $\chi_A \in A(G)$.
\begin{proposition}\label{qualrel}
Suppose that $G$ is a compact Abelian group and $A \subset G$ has
$\chi_A \in A(G)$. Then there is a finite $V \leq \wh{G}$ such
that $\chi_A=\chi_A \ast \mu_{V^\perp}$ almost everywhere.
\end{proposition}
This result is a special case of the idempotent theorem for
discrete groups first proved by Rudin in \cite{WRID}. The general
case (locally compact abelian groups) was proved by Cohen in
\cite{PC}.
\begin{proof}
[Proof of Proposition \ref{Cohen}] If $\chi_A \in A(G)$ then by
Proposition \ref{qualrel} there is some finite $V \leq \wh{G}$
such that $\chi_A = \chi_A \ast \mu_{V^\perp}$ almost everywhere.
For each $W \in G/V^\perp$ pick $x_W \in W$ such that $\chi_A \ast
\mu_{V^\perp}(x)=\chi_A(x_W)$ for almost all $x \in W$.
Integrating tells us that
\begin{equation*}
\alpha = \int{\chi_Ad\mu_G} = \int{\chi_A \ast
\mu_{V^\perp}d\mu_G}= |V|^{-1} \sum_{W \in G/V^\perp}{\chi_A(x_W)}
= n |V|^{-1}
\end{equation*}
for some integer $n$. It follows that $\{\alpha |V|\}(1-\{\alpha
|V|\}) = 0$ which contradicts the hypothesis on $\alpha$; hence
$\chi_A \not \in A(G)$.
\end{proof}

We use Bohr neighborhoods to prove Proposition \ref{qualrel};
these will also be useful in the proof of Theorem
\ref{maintheorem}.

\subsection{Bohr neighborhoods}\label{bohrsec} The following defines a natural valuation on $S^1$
\begin{equation*}
\|z\|:=\frac{1}{2\pi}\inf_{n \in \mathbb{Z}}{| 2\pi n + \arg z|},
\end{equation*}
which can be used to measure how far $\gamma(x)$ is from 1. This
leads to a definition of approximate annihilators for a finite
collection of characters $\Gamma$, namely
\begin{equation*}
B(\Gamma,\delta):=\{x \in G: \|\gamma(x)\| \leq \delta \textrm{
for all }\gamma \in \Gamma\};
\end{equation*}
such sets are called \emph{Bohr sets}, and their translates
\emph{Bohr neighborhoods}.

The following simple application of the pigeon-hole principle (the
details for which may be found, for example, in \cite{BJGSK})
shows that Bohr sets have positive measure.
\begin{lemma}\label{bohrsize}
Suppose $G$ is a compact Abelian group and $B(\Gamma,\delta)$ is a
Bohr set. Then $\mu_G(B(\Gamma,\delta)) \geq \delta^d$ where
$d:=|\Gamma|$.
\end{lemma}

Bohr sets are important because the translates of a Bohr sets
$B(\Gamma,\delta)$ are approximate joint level sets for the
characters in $\Gamma$ and hence a function that has the bulk of
its Fourier transform supported on a finite set $\Gamma$ is
approximately constant on all translates of $B(\Gamma,\delta)$.
Concretely we have the following lemma.
\begin{lemma}\label{contQual}
Suppose that $G$ is a compact Abelian group and $f \in A(G)$.
Suppose that $\eta \in (0,1]$. Then there is a finite set of
characters $\Gamma$ and a function $g$ such that
\begin{equation*}
\sup_{x-y \in B(\Gamma,\eta/3\|f\|_{A(G)})}{|g(x) - g(y)|} \leq
\eta \textrm{ and }f=g \textrm{ a.e.}
\end{equation*}
\end{lemma}
\begin{proof}
$f \in A(G)$ so there is a finite set of characters $\Gamma$ such
that
\begin{equation}\label{gamchoi}
\sum_{\gamma \not \in \Gamma}{|\wh{f}(\gamma)|} \leq \eta/3.
\end{equation}
Write
\begin{equation*}
g(x):=\sum_{\gamma \in \wh{G}}{\wh{f}(\gamma)\gamma(x)} \textrm{
and } g_\Gamma(x):=\sum_{\gamma \in
\Gamma}{\wh{f}(\gamma)\gamma(x)}.
\end{equation*}
Now
\begin{equation*}
|g_\Gamma(x) - g_\Gamma(y)| \leq \|f\|_{A(G)}\sup_{\gamma \in
\Gamma}{|\gamma(x)-\gamma(y)|} \leq \|f\|_{A(G)}\sup_{\gamma \in
\Gamma}{\|\gamma(x-y)\|} \leq \eta/3
\end{equation*}
if $x \in y+B(\Gamma,\eta/3\|f\|_{A(G)})$. From our choice of
$\Gamma$ in (\ref{gamchoi}) we have that $|g(x)-g_\Gamma(x)| \leq
\eta/3$ for all $x \in G$ so that
\begin{equation*}
|g(x) - g(y)| \leq |g(x) - g_\Gamma(x)| + |g_\Gamma(x) -
g_\Gamma(y)| + |g_\Gamma(y) - g(y)| \leq \eta
\end{equation*}
for all $x \in y+B(\Gamma,\eta/3\|f\|_{A(G)})$. By the inversion
formula $f=g$ almost everywhere, from which the lemma follows.
\end{proof}
If $f$ in the above lemma is the characteristic function of a set
then it takes only the values 0 or 1. We can use the following
version of the intermediate value theorem (implicit in
\cite{BJGSK} and \cite{TWSZP}) to show that such an $f$ must be
constant on cosets of $\langle B(\Gamma,\eta/3\|f\|_{A(G)})
\rangle$.
\begin{lemma}\label{IVT}\emph{(Discrete intermediate value theorem)}
Suppose that $G$ is a compact Abelian group and that
$B(\Gamma,\delta)$ is a Bohr set on $G$. Suppose that $g:G
\rightarrow \mathbb{R}$ has
\begin{equation}\label{conthyp}
\sup_{x-y \in B(\Gamma,\delta)}{|g(x) - g(y)|} \leq \eta.
\end{equation}
Suppose that $x_0,x_1 \in G$ have $x_0 - x_1 \in \langle
B(\Gamma,\delta) \rangle$. Then for any $c \in [g(x_0),g(x_1)]$,
there is an $x_2 \in G$ such that
\begin{equation*}
|g(x_2) - c| \leq \frac{\eta}{2}.
\end{equation*}
\end{lemma}
\begin{proof}
We write $H$ for the group, $\langle B(\Gamma,\delta) \rangle$,
generated by $B(\Gamma,\delta)$ and define
\begin{equation*}
S^- := \{x \in x_0 + H: g(x) < c - \frac{\eta}{2}\}
\end{equation*}
and
\begin{equation*}
S^+:=\{x \in x_0 + H: g(x) > c + \frac{\eta}{2}\}.
\end{equation*}
If the conclusion of the lemma is false then
$\mathcal{S}:=\{S^-,S^+\}$ is a partition of $x_0+H$.

By the continuity hypothesis (\ref{conthyp}) we have that if $x
\in S^-$ and $y \in B(\Gamma,\delta)$ then
\begin{equation*}
|g(x+y) - g(x)| \leq \eta \Rightarrow g(x+y) < c+\frac{\eta}{2}.
\end{equation*}
It follows that $x+y \not \in S^+$ and since $\mathcal{S}$ is a
partition of $x_0+H$ we conclude that $x+y \in S^-$. We have shown
that $S^- = S^- + H$.

Now $g(x_0) \leq c \leq g(x_1)$ and $\mathcal{S}$ is a partition
of $x_0+H$, whence $x_0 \in S^-$ and $x_1 \in S^+$. However $S^- =
S^- + H$, whence $S^- = x_0 + H$ and so $x_1 \in S^-$. This
contradicts the fact that $S^-$ and $S^+$ are disjoint and so
proves the lemma.
\end{proof}

\begin{proof}[Proof of Proposition \ref{qualrel}] Apply Lemma \ref{contQual} to $\chi_A$ with
$\eta=1/4$. The lemma provides a function $g$ with
\begin{equation}\label{conthyp2}
\sup_{x-y \in B(\Gamma,1/12\|\chi_A\|_{A(G)})}{|g(x) - g(y)|} \leq
1/4 \textrm{ and } g = \chi_A \textrm{ a.e.}
\end{equation}
Write $H:=\langle B(\Gamma,1/12\|\chi_A\|_{A(G)}) \rangle$.
\begin{claim}
$\chi_A$ is constant on cosets of $H$ (up to a null set).
\end{claim}
\begin{proof}
Suppose that $W$ is a coset of $H$ and there is a subset of $W$ of
positive measure on which $\chi_A$ is 0 and a subset of $W$ of
positive measure on which $\chi_A$ is 1. Since $g=\chi_A$ a.e. it
follows that there are $x_0,x_1 \in W$ (so that $x_0-x_1 \in H$)
with $g(x_0)=0$ and $g(x_1)=1$.

By Lemma \ref{IVT} there is some $x_2 \in G$ such that
$|g(x_2)-1/2| \leq 1/8$ and the fact that $g$ satisfies
(\ref{conthyp2}) then ensures that $|g(x) - 1/2| \leq 3/8$ for all
$x \in x_2 + B(\Gamma,1/12\|\chi_A\|_{A(G)})$. By Lemma
\ref{bohrsize} we know that $x_2 +
B(\Gamma,1/12\|\chi_A\|_{A(G)})$ as positive measure so that there
is some $x \in x_2 + B(\Gamma,1/12\|\chi_A\|_{A(G)})$ such that
$g(x)=\chi_A(x) \in \{0,1\}$. This contradicts the fact that $1/8
\leq g(x) \leq 7/8$. The claim follows.
\end{proof}
Again by Lemma \ref{bohrsize} we have that $\mu_G(H) \geq
\mu_G(B(\Gamma,1/12\|\chi_A\|_{A(G)}))>0$ so that $V:=H^\perp$ is
finite. Since $V^\perp=H$ the proposition is proved.
\end{proof}

The remainder of these notes establishes Theorem
\ref{maintheorem}. \S\ref{sec3} recalls some basic tools of local
Fourier analysis. Roughly \S\ref{finfourier} provides an effective
quantitative version of Lemma \ref{contQual} and
\S\ref{finphysical} establishes the necessary physical space
estimates to drive the result of \S\ref{finfourier}.
\S\ref{finrem} combines the results of the previous two sections
to prove the main result.

\section{Local Fourier analysis on compact Abelian
groups}\label{sec3}

Bourgain, in \cite{JB}, observed that one can localize the Fourier
transform to typical approximate level sets and retain approximate
versions of a number of the standard results for the Fourier
transform on compact Abelian groups. Since his original work
various expositions and extensions have appeared most notably in
the various papers of Green and Tao. We require a local version of
Chang's theorem as developed in \cite{TWSAS}; we follow the
preparatory discussion in there fairly closely.

\subsection{Approximate annihilators: typical Bohr sets and some of their
properties} We defined Bohr sets in \S\ref{bohrsec} and in view of
Lemma \ref{bohrsize} we write $\beta_{\Gamma,\delta}$, or simply
$\beta$ or $\beta_\delta$ if the parameters are implicit, for the
measure induced on $B(\Gamma,\delta)$ by $\mu_G$, normalised so
that $\|\beta_{\Gamma,\delta}\|=1$. Such measures are sometimes
referred to as \emph{normalised Bohr cutoffs}. We write $\beta'$
for $\beta_{\Gamma',\delta'}$, or $\beta_{\Gamma,\delta'}$ if no
$\Gamma'$ has been defined. Having defined these measures we norm
the $L^p$-spaces $L^p(B(\Gamma,\delta))$ in the obvious way
\emph{viz}.
\begin{equation*}
\|f\|_{L^p(B(\Gamma,\delta))}:=\left(\int{|f|^pd\beta_{\Gamma,\delta}}\right)^{\frac{1}{p}}.
\end{equation*}

As we noted Bohr sets can be thought of as approximate
annihilators, however genuine annihilators are also subgroups of
$G$, a property which, at least in an approximate form, we would
like to recover. Suppose that $\eta \in (0,1]$. Then
$B({\Gamma},\delta)+B({\Gamma},\eta\delta) \subset
B({\Gamma},(1+\eta)\delta)$. If $B({\Gamma},(1+\eta)\delta)$ is
not much bigger than $B({\Gamma},\delta)$ then we have a sort of
approximate additive closure in the sense that
$B({\Gamma},\delta)+B({\Gamma},\eta\delta) \approx
B({\Gamma},(1+\eta)\delta)$. Not all Bohr sets have this property
however Bourgain showed that typically they do. For our purposes
we have the following proposition.
\begin{proposition}\label{ubreg}
Suppose that $G$ is a compact Abelian group, $\Gamma$ a set of $d$
characters on $G$ and $\delta \in (0,1]$. There is an absolute
constant $c_{\mathcal{R}}>0$ and a $\delta' \in [\delta/2,\delta)$
such that
\begin{equation}\label{reg}
\frac{\mu_G(B(\Gamma,(1+\kappa)\delta'))}{\mu_G(B(\Gamma,\delta'))}
= 1 + O(|\kappa|d)
\end{equation}
whenever $|\kappa|d \leq c_{\mathcal{R}}$.
\end{proposition}
This result is not as easy as the rest of the section, it uses a
covering argument; a nice proof can be found in \cite{BJGTT}. We
say that $\delta'$ is \emph{regular for $\Gamma$} or that
$B(\Gamma,\delta')$ is \emph{regular} if
\begin{equation*}
\frac{\mu_G(B(\Gamma,(1+\kappa)\delta'))}{\mu_G(B(\Gamma,\delta'))}
= 1 + O(|\kappa|d) \textrm{ whenever } |\kappa|d \leq
c_{\mathcal{R}}.
\end{equation*}
It is regular Bohr sets to which we localize the Fourier transform
and we begin by observing that normalised regular Bohr cutoffs are
approximately translation invariant and so function as normalised
approximate Haar measures.
\begin{lemma}
\label{contlem}\emph{(Normalized approximate Haar measure)}
Suppose that $G$ is a compact Abelian group and $B(\Gamma,\delta)$
is a regular Bohr set. If $y \in B(\Gamma,\delta')$ then
$\|(y+\beta_\delta) - \beta_\delta\| \ll d\delta'\delta^{-1}$
where we recall that $y+\beta_\delta$ denotes the measure
$\beta_\delta$ composed with translation by $y$.
\end{lemma}
The proof follows immediately from the definition of regularity.
In applications the following corollary will be useful but it
should be ignored until it is used.
\begin{corollary}
\label{contlm} Suppose that $G$ is a compact Abelian group and
$B(\Gamma,\delta)$ is a regular Bohr set. If $f \in L^\infty(G)$
then
\begin{equation*}
\|f\ast \beta - f \ast \beta(x)\|_{L^\infty(x+B(\Gamma,\delta'))}
\ll \|f\|_{L^\infty(G)}d\delta'\delta^{-1}.
\end{equation*}
\end{corollary}

With an approximate Haar measure we are in a position to define
the local Fourier transform: Suppose that $x'+B(\Gamma,\delta)$ is
a regular Bohr neighborhood. Then we define the Fourier transform
local to $x'+B(\Gamma,\delta)$ by
\begin{equation*}
L^1(x'+B(\Gamma,\delta)) \rightarrow \ell^\infty(\wh{G}); f
\mapsto \wh{fd(x'+\beta_{\Gamma,\delta})}.
\end{equation*}

\subsection{The structure of sets of characters supporting large
values of the local Fourier transform} Having defined the local
transform and recorded the key tools it remains for us to recall
the result from \cite{TWSAS} to which the title of these notes
refers.
\begin{proposition}\label{local Changs bound} \emph{(Chang's theorem local to Bohr sets, \cite{TWSAS}, Proposition 5.2)}
Suppose that $G$ is a compact Abelian group and $B(\Gamma,\delta)$
is a regular Bohr set. Suppose that $f \in L^2(B(\Gamma,\delta))$
and $\epsilon,\eta \in (0,1]$. Then there is a set of characters
$\Lambda$ and a $\delta' \in (0,1]$ with
\begin{equation*}
|\Lambda| \ll \epsilon^{-2}\log
\|f\|_{L^1(B(\Gamma,\delta))}^{-2}\|f\|_{L^2(B(\Gamma,\delta))}^2
\end{equation*}
and
\begin{equation*}
\delta' \gg \delta\eta\epsilon^2/d^2 \log
\|f\|_{L^1(B(\Gamma,\delta))}^{-2}\|f\|_{L^2(B(\Gamma,\delta))}^2
\end{equation*}
and furthermore,
\begin{equation*}
\{\gamma \in \wh{G}:|\wh{fd\beta}(\gamma)| \geq \epsilon
\|f\|_{L^1(B(\Gamma,\delta))}\}
\end{equation*}
is contained in
\begin{equation*}
\{\gamma \in \wh{G}:|1-\gamma(x)| \leq \eta \textrm{ for all } x
\in B(\Gamma \cup \Lambda,\delta')\}.
\end{equation*}
\end{proposition}

\section{An iteration argument in Fourier space}\label{finfourier}

The main result of this section takes physical space information
about a set $A \subset G$ and converts it into Fourier
information. The lemma is based on Lemma 4.8 in \cite{TWSCVS} with
two main modifications:
\begin{itemize}
\item We have to assume the comparability of the local $L^2$-norm
squared and local $L^1$-norm; ensuring this hypothesis is the
principal extra complication of \S\ref{finphysical}. \item We are
less careful in our analysis because the physical space estimates
available to us in the general setting are sufficiently weak as to
render any more care irrelevant.
\end{itemize}

\begin{lemma}\label{finitlem}
\emph{(Iteration lemma)} Suppose that $G$ is a compact Abelian
group, $B(\Gamma,\delta)$ is a Bohr set and $B(\Gamma,\delta')$ is
a regular Bohr set. Suppose that $A \subset G$ has $\chi_A \in
A(G)$ and write $f:=\chi_A - \chi_A \ast \beta$. Suppose,
additionally, that
\begin{equation*}
\|f\|_{L^2(B(\Gamma,\delta'))}^2 \asymp
\|f\|_{L^1(B(\Gamma,\delta'))} \textrm{ and }
\|f\|_{L^2(B(\Gamma,\delta'))}^2 >0.
\end{equation*}
Suppose that $\epsilon \in (0,1]$ is a parameter. Then either
$\|\chi_A\|_{A(G)} \gg \epsilon^{-1}$ or there is a set of
characters $\Lambda$ and a regular Bohr set
$B(\Gamma\cup\Lambda,\delta'')$ such that
\begin{equation*}
|\Lambda| \ll \epsilon^{-2}\log
\|f\|_{L^2(B(\Gamma,\delta'))}^{-1} \textrm{ and } \delta'' \gg
\delta'\epsilon^3 /d^2 \log \|f\|_{L^2(B(\Gamma,\delta'))}^{-1},
\end{equation*}
where, as usual, $d:=|\Gamma|$, and
\begin{equation*}
\sum_{\gamma \in \mathcal{N} \setminus
\mathcal{O}}{|\wh{\chi_A}(\gamma)|} \gg 1,
\end{equation*}
where $\mathcal{O}:=\{\gamma:|1-\gamma(x)| \leq \epsilon \textrm{
for all } x \in B(\Gamma,\delta)\}$ and
$\mathcal{N}:=\{\gamma:|1-\gamma(x)| \leq \epsilon \textrm{ for
all } x \in B(\Gamma\cup \Lambda,\delta'')\}$.
\end{lemma}
\begin{proof}
By Plancherel's theorem we have
\begin{equation}\label{planch}
\sum_{\gamma
\in\wh{G}}{\wh{f}(\gamma)\overline{\wh{fd\beta'}(\gamma)}} =
\|f\|_{L^2(B(\Gamma,\delta'))}^2.
\end{equation}
Write
\begin{equation*}
\mathcal{L}:=\{\gamma: |\wh{fd\beta'}(\gamma)| \geq \epsilon
\|f\|_{L^1(B(\Gamma,\delta'))}\},
\end{equation*}
and suppose that
\begin{equation}\label{hyp}
\sum_{\gamma \not \in
\mathcal{L}}{\wh{f}(\gamma)\overline{\wh{fd\beta'}(\gamma)}} \geq
\|f\|_{L^2(B(\Gamma,\delta'))}^2/2.
\end{equation}
Note that
\begin{equation*}
\|f\|_{A(G)} = \|\chi_A - \chi_A \ast \beta\|_{A(G)} \leq
\|\chi_A\|_{A(G)} + \|\chi_A \ast \beta\|_{A(G)} \leq
2\|\chi_A\|_{A(G)},
\end{equation*}
whence
\begin{eqnarray*}
\sum_{\gamma \not \in
\mathcal{L}}{|\wh{f}(\gamma)||\wh{fd\beta'}(\gamma)|} & \leq &
\epsilon \|f\|_{L^1(B(\Gamma,\delta'))} \|f\|_{A(G)}\\ & \leq &
2\epsilon\|\chi_A\|_{A(G)} \|f\|_{L^1(B(\Gamma,\delta'))}\\ & \ll
& \epsilon \|\chi_A\|_{A(G)} \|f\|_{L^2(B(\Gamma,\delta'))}^2.
\end{eqnarray*}
If (\ref{hyp}) holds then the left hand side of this is at least
$\|f\|_{L^2(B(\Gamma,\delta'))}^2/2$ and so (dividing by
$\|f\|_{L^2(B(\Gamma,\delta'))}^2$) we conclude that
$\|\chi_A\|_{A(G)} \gg \epsilon^{-1}$.

Thus we may suppose that (\ref{hyp}) is not true and therefore, by
(\ref{planch}), that
\begin{equation*}
\sum_{\gamma  \in
\mathcal{L}}{\wh{f}(\gamma)\overline{\wh{fd\beta'}(\gamma)}} \geq
\|f\|_{L^2(B(\Gamma,\delta'))}^2/2.
\end{equation*}
By Proposition \ref{local Changs bound} there is a set of
characters $\Lambda$ and a $\delta''$ (regular for $\Gamma \cup
\Lambda$ by Proposition \ref{ubreg}), with
\begin{equation*}
|\Lambda| \ll \epsilon^{-2}\log
\|f\|_{L^1(B(\Gamma,\delta'))}^{-2}\|f\|_{L^2(B(\Gamma,\delta'))}^2
\ll \epsilon^{-2} \log \|f\|_{L^2(B(\Gamma,\delta'))}^{-1}
\end{equation*}
and
\begin{equation*}
\delta'' \gg \delta'\epsilon^3/d^2 \log
\|f\|_{L^1(B(\Gamma,\delta'))}^{-2}\|f\|_{L^2(B(\Gamma,\delta'))}^2
\gg \delta'\epsilon^3/d^2\log \|f\|_{L^2(B(\Gamma,\delta'))}^{-1},
\end{equation*}
such that
\begin{equation*}
\mathcal{L} \subset \{\gamma:|1-\gamma(x)| \leq \epsilon \textrm{
for all } x \in B(\Gamma\cup\Lambda,\delta'')\}=\mathcal{N}.
\end{equation*}
Since $\mathcal{L} \subset \mathcal{N}$ we have
\begin{equation*}
\sum_{\gamma \in
\mathcal{N}}{|\wh{f}(\gamma)\wh{fd\beta'}(\gamma)|} \geq
\|f\|_{L^2(B(\Gamma,\delta'))}^2/2.
\end{equation*}
Now
\begin{equation*}
|\wh{fd\beta'}(\gamma)| \leq \|f\|_{L^1(B(\Gamma,\delta'))} \ll
\|f\|_{L^2(B(\Gamma,\delta'))}^2,
\end{equation*}
hence
\begin{equation}\label{Nlarge}
\sum_{\gamma \in \mathcal{N}}{|\wh{f}(\gamma)|} \gg 1.
\end{equation}
Finally suppose that
\begin{equation}\label{bigO}
\sum_{\gamma \in \mathcal{O}}{|\wh{f}(\gamma)|} \geq
\frac{1}{2}\sum_{\gamma \in \mathcal{N}}{|\wh{f}(\gamma)|}.
\end{equation}
By the definition of $\mathcal{O}$ we have
\begin{equation*}
\sum_{\gamma \in \mathcal{O}}{|\wh{f}(\gamma)|} = \sum_{\gamma \in
\mathcal{O}}{|\wh{\chi_A}(\gamma)|.|1-\wh{\beta}(\gamma)|} \leq
\|\chi_A\|_{A(G)}\sup_{\gamma \in
\mathcal{O}}{|1-\wh{\beta}(\gamma)|} \leq \epsilon
\|\chi_A\|_{A(G)}.
\end{equation*}
It follows that if (\ref{bigO}) holds then, in view of
(\ref{Nlarge}), $\|\chi_A\|_{A(G)} \gg \epsilon^{-1}$. Thus we may
assume it does not and hence that
\begin{equation*}
\sum_{\gamma \in \mathcal{N} \setminus
\mathcal{O}}{|\wh{f}(\gamma)|} \gg 1.
\end{equation*}
Noting that $|\wh{f}(\gamma)| \leq 2 |\wh{\chi_A}(\gamma)|$
completes the proof.
\end{proof}

\section{Physical space estimates}\label{finphysical}

The objective of this section is to prove the following result.
\begin{proposition}\label{sizeprop} Suppose that $G$ is a finite Abelian
group and $B(\Gamma,\delta)$ is a regular Bohr set in $G$. Suppose
that $A \subset G$ has density $\alpha$ and for all finite $V \leq
\wh{G}$ with $|V| \leq M$ we have
$\{\alpha|V|\}(1-\{\alpha|V|\})\gg 1$. Then either
\begin{equation*}
\log M \ll d(\log \delta^{-1}+d\log d)
\end{equation*}
or there is an $x'' \in G$ and reals $\delta'$ and $\delta''$,
both regular for $\Gamma$, with $\delta' \leq \delta$,
\begin{equation*}
\log \delta \delta'^{-1} \ll d \log d \textrm{ and } \log
\delta''^{-1} \ll d (\log \delta^{-1} + d \log d)
\end{equation*}
such that
\begin{equation}\label{setlike}
\|\chi_A - \chi_A \ast \beta'\|_{L^2(x''+B(\Gamma,\delta''))}^2
\asymp \|\chi_A - \chi_A \ast
\beta'\|_{L^1(x''+B(\Gamma,\delta''))}
\end{equation}
and
\begin{equation}\label{lwrbd}
\log \|\chi_A - \chi_A \ast
\beta'\|_{L^2(x''+B(\Gamma,\delta''))}^{-2} \ll d(\log \delta^{-1}
+ d \log d).
\end{equation}
\end{proposition}
Of the two parts (\ref{setlike}) and (\ref{lwrbd}) the second is
the easiest to derive and comes essentially from a straightforward
generalization of the physical space estimates of \cite{TWSCVS}
combined with discrete intermediate value theorem (Lemma
\ref{IVT}). To ensure (\ref{setlike}) requires more work and is
the principal extra ingredient of these notes.

We begin with a version of Lemma 5.1, \cite{TWSCVS}, appropriate
to our more general setting. In fact the proof which follows is
slightly simpler than that in \cite{TWSCVS} and would have been
sufficient for application there as well; the weakness of the
present approach only impacts on the implied constants.
\begin{lemma}\label{bse}
Suppose that $G$ is a finite Abelian group. Suppose that $f\in
L^1(G)$ maps $G$ into $[0,1]$ and that $V \leq \wh{G}$, finite,
has $\{\|f\|_1 |V|\}(1-\{\|f\|_1|V|\}) \gg 1$. Then there is a
coset $x'+V^\perp$ with
\begin{equation*}
f \ast \mu_{V^\perp}(x') \gg \mu_G(V^\perp) \textrm{ and } (1-f)
\ast \mu_{V^\perp}(x') \gg \mu_G(V^\perp).
\end{equation*}
\end{lemma}
\begin{proof}
$f \ast \mu_{V^\perp}$ is constant on cosets of $V^\perp$ so we
define
\begin{equation*}
g(x):=\begin{cases} 1 & \textrm{ if } f \ast \mu_{V^\perp}(x) \geq
1/2\\ 0 & \textrm{ otherwise.}
\end{cases}
\end{equation*}
Since $g$ is integral on cosets of $V^\perp$ there is some integer
$n$ such that
\begin{equation*}
\int{gd\mu_G} = n\mu_G(V^\perp).
\end{equation*}
However
\begin{eqnarray*}
|n\mu_G(V^\perp) - \|f\|_1| & = & |\int{gd\mu_G} - \int{fd\mu_G}|\\
& = & |\int{gd\mu_G} - \int{f \ast \mu_{V^\perp}d\mu_G}|\\ & \leq
& \int{|g - f \ast \mu_{V^\perp}|d\mu_G}\\ & \leq &\sup_{x \in
G}{\min\{f \ast \mu_{V^\perp}(x),1- f \ast \mu_{V^\perp}(x)\}}\\ &
= & \sup_{x \in G}{\min\{f \ast \mu_{V^\perp}(x),(1- f )\ast
\mu_{V^\perp}(x)\}}.
\end{eqnarray*}
Now
\begin{eqnarray*}
|n\mu_G(V^\perp) - \|f\|_1| & = & \mu_G(V^\perp)|n-|V|\|f\|_1|\\ &
\geq & \mu_G(V^\perp)\{\|f\|_1|V|\}(1-\{\|f\|_1|V|\})\\ & \gg &
\mu_G(V^\perp),
\end{eqnarray*}
and the conclusion of the lemma follows.
\end{proof}
We require one more preliminary lemma.
\begin{lemma}\label{gplemma}
Suppose that $G$ is a finite Abelian group and $B(\Gamma,\delta)$
is a Bohr set in $G$. Then there are reals $\delta'$ and
$\delta''$ both regular for $\Gamma$ with $\delta' \leq \delta$,
\begin{equation*}
\log \delta\delta'^{-1} \ll d \log d \textrm{ and } \delta'' \gg
\delta'/d
\end{equation*}
such that
\begin{equation*}
B(\Gamma,\delta')^\perp = B(\Gamma,\delta'')^\perp
\end{equation*}
and
\begin{equation*}
\|f\ast \beta' - f \ast
\beta'(x)\|_{L^\infty(x+B(\Gamma,\delta''))} \leq \|f\|_\infty/4
\end{equation*}
for all $x \in G$ and $f \in L^\infty(G)$.
\end{lemma}
\begin{proof} We define a sequence
$(\delta_i)_i$ iteratively and write
\begin{equation*}
\beta_i:=\beta_{\Gamma,\delta_i} \textrm{ and } H_i:=\langle
B(\Gamma,\delta_i)\rangle.
\end{equation*}
To begin with we apply Proposition \ref{ubreg} to get some
$\delta_0$ regular for $\Gamma$ with $\delta \geq \delta_0 \gg
\delta$. Now, if we have constructed $\delta_i$ for some $i \geq
0$, we apply Corollary \ref{contlm} (and Proposition \ref{ubreg})
to get a $\delta_{i+1}$ regular for $\Gamma$ with $\delta_i \geq
\delta_{i+1} \gg \delta_i/d$ and
\begin{equation*}
\|f \ast \beta_i - f \ast
\beta_i(x)\|_{L^\infty(x+B(\Gamma,\delta_{i+1}))} \leq
\|f\|_\infty/4
\end{equation*}
for all $x \in G$ and $f \in L^\infty(G)$. We are done if we can
show that there is some $i \leq d$ such that $H_i=H_{i+1}$. This
follows by the pigeon-hole principle from the following claim.

\begin{claim} Suppose that $\kappa_0 \in (0,1]$. Then there is a
sequence of elements $x_1,...,x_d \in G$ such that for each
$\kappa \in (\kappa_0,1]$ there is some $0 \leq i \leq d$ such
that
\begin{equation*}
\langle B(\Gamma,\kappa)\rangle =\langle x_1,...,x_i \rangle +
\bigcap_{\gamma \in \Gamma}{\ker \gamma}.
\end{equation*}
\end{claim}
\begin{proof}
The proof of the claim is based on ideas from the geometry of
numbers introduced to the area by Ruzsa in \cite{IZR};
\cite{BJGIZR} contains a neat exposition of the idea. By
quotienting we may assume that $\bigcap_{\gamma \in \Gamma}{\ker
\gamma} = \{0_G\}$.

Let $\phi:G \rightarrow \mathbb{T}^d; x \mapsto
(\gamma(x))_{\gamma \in \Gamma}$ and define the lattice
$\mathcal{L}:=\bigcup{\phi(G)} \leq \mathbb{R}^d$. Since
$\bigcap_{\gamma \in \Gamma}{\ker \gamma} = \{0_G\}$ there is a
natural homomorphism $\psi:\mathcal{L} \rightarrow G$ which takes
$b \in \mathcal{L}$ to the unique $x \in G$ such that
$\phi(x)=b+\mathbb{Z}^d$, with kernel $\mathbb{Z}^d$.

We write $Q$ for the unit cube centered at the origin in
$\mathbb{R}^d$ and note that $\psi(\kappa Q)= B(\Gamma,\kappa)$.
We choose linearly independent vectors $b_1,..,b_d \in
\mathcal{L}$ inductively so that
\begin{equation*}
\|b_i\|_\infty \leq \inf\{\lambda:\lambda Q \cap \mathcal{L}
\textrm{ contains } i \textrm{ linearly independent vectors}\}.
\end{equation*}
Let $x_i = \psi(b_i)$. Since $\psi$ is linear we have $\langle
B(\Gamma,\kappa)\rangle = \psi(\langle \kappa Q \rangle)$, but to
each $\kappa \in (0,1]$ there corresponds an $1 \leq i \leq d$
such that $\kappa Q$ contains at most $i$ linearly independent
vectors and $\kappa Q$ contains $b_1,...,b_i$. Hence $\langle
\kappa Q \rangle = \langle b_1,...,b_i \rangle$. Linearity of
$\psi$ again gives the result.
\end{proof}
\end{proof}

\begin{proof}[Proof of Proposition \ref{sizeprop}.]
Applying Lemma \ref{gplemma} we get reals $\delta'\leq \delta$ and
$\delta'''$ both regular for $\Gamma$ with
\begin{equation*}
\log \delta\delta'^{-1} \ll d\log d \textrm{ and } \delta''' \gg
\delta'/d
\end{equation*}
such that
\begin{equation*}
B(\Gamma,\delta')^\perp = B(\Gamma,\delta''')^\perp =: V,
\end{equation*}
and
\begin{equation}\label{contfct}
\|\chi_A\ast \beta' - \chi_A \ast
\beta'(x)\|_{L^\infty(x+B(\Gamma,\delta'''))} \leq 1/4 \textrm{
for all } x \in G.
\end{equation}

Now
\begin{equation}\label{tgy}
|V|^{-1} = \mu_G(V^\perp) \geq \mu_G(B(\Gamma,\delta')) \geq
(\delta')^d;
\end{equation}
the last inequality by Lemma \ref{bohrsize}. If $M \leq |V|$ then
we are in the first case of the lemma. Otherwise by hypothesis
$\{\alpha |V|\}(1-\{\alpha |V|\}) \gg 1$. If we put $f=\chi_A \ast
\beta'$ then $\|f\|_1 = \alpha$ and $f$ maps $G$ into $[0,1]$
whence, by Lemma \ref{bse}, there is some $x''' \in G$ such that
\begin{equation}\label{bhg}
\chi_A \ast \beta'\ast \mu_{V^\perp}(x''') \gg \mu_G(V^\perp)
\textrm{ and } (1-\chi_A \ast \beta')\ast \mu_{V^\perp}(x''') \gg
\mu_G(V^\perp).
\end{equation}

The argument now splits into three cases.
\begin{enumerate}
\item \emph{There are elements $x_0,x_1 \in x'''+V^\perp$ such
that $\chi_A \ast \beta'(x_0) \geq 1/2$ and $\chi_A \ast
\beta'(x_1) \leq 1/2$}. Here $x_0 - x_1 \in V^\perp =
B(\Gamma,\delta''')^{\perp\perp} = \langle B(\Gamma,\delta''')
\rangle$, so by the discrete intermediate value theorem (Lemma
\ref{IVT}) and (\ref{contfct}) we conclude that there is some $x_2
\in x'''+ V^\perp$ such that
\begin{equation*}
\frac{3}{8} \leq \chi_A \ast \beta'(x_2) \leq \frac{5}{8}.
\end{equation*}
Further by (\ref{contfct}) we conclude that
\begin{equation*}
\frac{1}{8} \leq \chi_A \ast \beta'(x) \leq \frac{7}{8} \textrm{
for all } x \in x_2 + B(\Gamma,\delta''').
\end{equation*}
Since $\chi_A$ only takes values in $\{0,1\}$ it follows that
$|\chi_A - \chi_A \ast \beta'| \asymp 1$ on
$x_2+B(\Gamma,\delta''')$. Thus
\begin{equation*}
\|\chi_A - \chi_A\ast \beta'\|_{L^2(x_2+B(\Gamma,\delta'''))}^2
\asymp \mu_G(B(\Gamma,\delta''')),
\end{equation*}
and
\begin{equation*}
\|\chi_A - \chi_A\ast \beta'\|_{L^1(x_2+B(\Gamma,\delta'''))}
\asymp \mu_G(B(\Gamma,\delta''')).
\end{equation*}
The result follows on putting $\delta''=\delta'''$; Lemma
\ref{bohrsize} then gives (\ref{lwrbd}).

\item \emph{$\chi_A \ast \beta'(x) \leq 1/2$ for all $x \in x'''+
V^\perp$}. Suppose that $\delta'' \leq \delta'$. Then
$B(\Gamma,\delta'') \subset B(\Gamma,\delta')$ so we have $\langle
B(\Gamma,\delta'') \rangle \subset \langle B(\Gamma,\delta')
\rangle$, whence $\beta'' \ast \mu_{V^\perp} = \mu_{V^\perp}=
\beta' \ast \mu_{V^\perp}$. Thus we define
\begin{equation*}
\alpha':=\chi_A \ast \beta' \ast \mu_{V^\perp}(x''') = \chi_A \ast
\beta''\ast \mu_{V^\perp}(x'''),
\end{equation*}
which has
\begin{equation*}
\alpha' \gg \mu_G(V^\perp) \geq (\delta')^d)^d>0
\end{equation*} by (\ref{bhg}) and (\ref{tgy}).
By Corollary \ref{contlm} we can pick a $\delta''$ (regular for
$\Gamma$ by Proposition \ref{ubreg}) with $\delta' \geq \delta''
\gg \delta'\alpha'/d$ such that
\begin{equation}\label{contdl}
\|\chi_A - \chi_A \ast
\beta'(x)\|_{L^\infty(x+B(\Gamma,\delta''))} \leq \alpha' \textrm{
for all } x \in G.
\end{equation}

We write
\begin{equation*}
L:=\{x \in x''' + V^\perp: \chi_A \ast \beta''(x) \geq \alpha'/2\}
\end{equation*}
and note that
\begin{equation*}
\int_{x \not \in L}{\chi_A \ast \beta''(x)d\mu_{V^\perp}(x'''-x)}
\leq \sup_{x \not \in L}{\chi_A \ast \beta''(x)} \leq \alpha'/2,
\end{equation*}
so
\begin{equation*}
\int_{x \in L}{\chi_A \ast \beta''(x)d\mu_{V^\perp}(x'''-x)} \geq
\alpha'/2.
\end{equation*}
If $\chi_A \ast \beta''(x) \neq 0$ then $\chi_A \ast \beta'(x)
\neq 0$, whence
\begin{eqnarray*}
\alpha'/2 & \leq & \int_{x \in L}{\chi_A \ast
\beta'(x)\frac{\chi_A \ast \beta''(x)}{\chi_A \ast
\beta'(x)}d\mu_{V^\perp}(x'''-x)}\\ & \leq & \alpha' \sup_{x \in
L}{\frac{\chi_A \ast \beta''(x)}{\chi_A \ast \beta'(x)}},
\end{eqnarray*}
Dividing by $\alpha'$ (which we have previously observed is
positive) we conclude that there is some $x'' \in L$ such that
\begin{equation*}
\chi_A \ast \beta''(x'') \geq \chi_A \ast \beta'(x'')/4.
\end{equation*}
If $x \in A \cap (x''+B(\Gamma,\delta''))$ then $|\chi_A(x) -
\chi_A \ast \beta'(x)| \asymp 1$ since $\chi_A \ast \beta'(x) \leq
1/2$ by the hypothesis of this case. If $x \in A^c \cap (x'' +
B(\Gamma,\delta''))$ then
\begin{equation*}
|\chi_A(x) - \chi_A \ast \beta'(x)| \leq |\chi_A \ast \beta'(x)|
\leq \chi_A \ast \beta'(x'') + O(\alpha') = O(\alpha'),
\end{equation*}
where the second inequality is a result of (\ref{contdl}). It
follows that
\begin{equation*}
\|\chi_A - \chi_A \ast \beta'\|_{L^1(x'' + B(\Gamma,\delta''))}
\gg \chi_A \ast \beta''(x''),
\end{equation*}
and
\begin{eqnarray*}
\|\chi_A - \chi_A \ast \beta'\|_{L^1(x'' + B(\Gamma,\delta''))} &
\leq & O(\chi_A \ast \beta''(x'') ) + O(\alpha')\\ & = & O(\chi_A
\ast \beta''(x''))
\end{eqnarray*}
since $\chi_A \ast \beta''(x'') \gg \alpha'$ since $x'' \in L$.
Similarly we have
\begin{equation*}
\|\chi_A - \chi_A \ast \beta'\|_{L^2(x'' + B(\Gamma,\delta''))}^2
\gg \chi_A \ast \beta''(x''),
\end{equation*}
and
\begin{eqnarray*}
\|\chi_A - \chi_A \ast \beta'\|_{L^2(x'' + B(\Gamma,\delta''))}^2
& \leq & O(\chi_A \ast \beta''(x'') ) + O(\alpha')\\ & = &
O(\chi_A \ast \beta''(x'')).
\end{eqnarray*}
It follows that
\begin{equation*}
\|\chi_A - \chi_A \ast \beta'\|_{L^1(x'' + B(\Gamma,\delta''))}
\asymp \|\chi_A - \chi_A \ast \beta'\|_{L^2(x'' +
B(\Gamma,\delta''))}^2,
\end{equation*}
and
\begin{equation*}
\|\chi_A - \chi_A \ast \beta'\|_{L^2(x'' + B(\Gamma,\delta''))}^2
\gg \alpha'.
\end{equation*}

\item \emph{$\chi_A \ast \beta'(x) \geq 1/2$ for all $x \in
x'''+V^\perp$}. This follows by replacing $A$ in the previous case
by $A^c$.
\end{enumerate}
The proof is complete.
\end{proof}

\section{Proof of Theorem \ref{maintheorem}}\label{finrem}

\begin{proof} In what follows it is convenient to let $C>0$ denote
an absolute constant which may vary from instance to instance.

Fix $\epsilon \in (0,1]$ to be optimized later. We define three
sequences $(\delta_k)_k$, $(\delta_k')_k$ and $(\delta_k'')_k$ of
reals, one sequence $(x_k)_k$ of elements of $G$, and one sequence
$(\Gamma_k)_k$ of sets of characters inductively. We write
\begin{equation*}
\beta_k:=\beta_{\Gamma_k,\delta_k},
\beta_k':=\beta_{\Gamma_k,\delta_k'} \textrm{ and }
\beta_k'':=\beta_{\Gamma_k,\delta_k''},
\end{equation*}
as well as $d_k:=|\Gamma_k|$ and
\begin{equation*}
\mathcal{L}_k:=\{\gamma:|1-\gamma(x)| \leq \epsilon \textrm{ for
all }x \in B(\Gamma_k,\delta_k)\}.
\end{equation*}
We shall ensure the following properties.
\begin{enumerate}
\item \label{p1} $B(\Gamma_k,\delta_k)$, $B(\Gamma_k,\delta_k')$
and $B(\Gamma_k,\delta_k'')$ are regular; \item \label{p2}
\begin{equation*} \|\chi_A - \chi_A\ast
\beta_k\|_{L^2(x_k+B(\Gamma_k,\delta_k'))}^2 \asymp \|\chi_A -
\chi_A \ast \beta_k'\|_{L^1(x_k+B(\Gamma_k,\delta_k''))};
\end{equation*} \item \label{p3} \begin{equation*}
\|\chi_A - \chi_A \ast
\beta_k'\|_{L^2(x_k+B(\Gamma_k,\delta_k''))}^2 \gg
\delta_k^{d_k}/(Cd_k)^{d_k^2};
\end{equation*}
\item \label{p4} \begin{equation*} \delta_k \geq \delta_k' \gg
\delta_k/(Cd_k)^{d_k} \textrm{ and } \delta_k \geq \delta_k'' \gg
\delta_k^{d_k+1}/(Cd_k)^{d_k^2+d_k+1};
\end{equation*}
\item \label{p5}
\begin{equation*}
d_{k+1} \ll \epsilon^{-2}d_k (\log \delta_k^{-1} + d_k \log d_k);
\end{equation*}
\item \label{p6}
\begin{equation*}
\delta_k \geq \delta_{k+1} \gg
\delta_k^{d_k+1}\epsilon^4/(Cd_k)^{d_k^2+d_k+6}\log \delta_k^{-1};
\end{equation*}
\item \label{p7} \begin{equation*} \sum_{\gamma \in
\mathcal{L}_{k+1}\setminus \mathcal{L}_k}{|\wh{\chi_A}(\gamma)|}
\gg 1.\end{equation*}
\end{enumerate}

We initialize the iteration with $\Gamma_0=\{0_{\wh{G}}\}$. Pick
$\delta_0 \gg 1$ regular for $\Gamma_0$ by Proposition
\ref{ubreg}. Apply Proposition \ref{sizeprop} (assuming that
$\delta_0^{d_0}(Cd_0)^{d_0^2} < M$) to get $x_0$ $\delta_0'$ and
$\delta_0''$ satisfying properties
(\ref{p1}),(\ref{p2}),(\ref{p3}) and (\ref{p4}). By translating
$A$ by $-x_0$, if necessary, we can apply Lemma \ref{finitlem}
(assuming that we don't have $\|\chi_A\|_{A(G)} \gg
\epsilon^{-1}$) to get $\Gamma_1$ and $\delta_1$ such that
properties (\ref{p5}), (\ref{p6}) and (\ref{p7}) are satisfied.

Given $\Gamma_k$ and $\delta_k$ we can proceed as we just have
(assuming that $\delta_k^{d_k}(Cd_k)^{d_k^2} < M$) to generate
$x_k$, $\delta_k'$, $\delta_k''$, $\delta_{k+1}$ and
$\Gamma_{k+1}$.

By property (\ref{p7}) (and the leftmost inequality in (\ref{p6}))
we have $\|\chi_A\|_{A(G)} \gg k$, so either $\|\chi_A\|_{A(G)}
\gg \epsilon^{-1}$, or the iteration terminates with $k \ll
\epsilon^{-1}$.

(\ref{p5}) and (\ref{p6}) imply
\begin{equation*}
d_{k+1} \ll \epsilon^{-2} d_k \log \delta_k^{-1} \textrm{ and }
\log \delta_{k+1}^{-1} \ll d_k \log \delta_k^{-1} + \log
\epsilon^{-1},
\end{equation*}
whence
\begin{equation*}
d_{k+1} \ll \epsilon^{-2} d_k\log \delta_k^{-1} \textrm{ and }
\epsilon^{-2} \log \delta_{k+1}^{-1} \ll \epsilon^{-2} d_k \log
\delta_k^{-1}.
\end{equation*}
It follows that $d_{k+1} \ll d_k^2$ and so $d_k \leq 2^{2^{Ck}}$
and $\delta_k \geq 2^{2^{2^{Ck}}}$. For the iteration to terminate
we must have $M \leq \delta_k^{d_k}(Cd_k)^{d_k^2}$; for this to
happen for some $k \ll \epsilon^{-1}$ we need
$2^{2^{2^{C\epsilon^{-1}}}} \geq M$. The result follows.
\end{proof}

\section{Acknowledgements} I should like to thank Tim Gowers and Ben Green for supervision
and much encouragement.


\begin{thebibliography}{9}% Replace 9 by 99 if 10 or more references
\bibitem{JB}
J.~Bourgain, \emph{On triples in arithmetic progression}. GAFA 9
(1999), no. 5, 968--984.
\bibitem{JB2}
J.~Bourgain, \emph{On the distributions of the Fourier spectrum of
Boolean functions}. Israel J. Math. 131 (2002), 269--276.
\bibitem{MCC}
M-C.~Chang, \emph{A polynomial bound in Freiman's theorem}. Duke
Math. J. 113 (2002), no. 3, 399--419.
\bibitem{PC}
Cohen
\bibitem{BJGSK}
B.J.~Green and S.~Konyagin, \emph{On the Littlewood problem modulo
a prime}. Preprint.
\bibitem{BJGIZR}
B.J.~Green and I.Z.~Ruzsa, \emph{Freiman's Theorem in an arbitrary
abelian group}. Preprint.
\bibitem{BJGTT}
B.J.~Green and T.~Tao, \emph{An Inverse Theorem for the Gowers
$U^3$-norm}. Preprint.
\bibitem{WRID}
W.~Rudin, \emph{Idempotent measures on Abelian groups}. Pacific J.
Math. 9 1959 195--209.
\bibitem{WR}
W.~Rudin, \emph{Fourier analysis on groups}. Interscience Tracts
in Pure and Applied Mathematics, No. 12 Interscience Publishers (a
division of John Wiley and Sons), New York-London 1962 ix+285 pp.
\bibitem{IZR}
I.Z.~Ruzsa, \emph{Generalized arithmetical progressions and
sumsets}. Acta Math. Hungar. 65 (1994), no. 4, 379--388.
\bibitem{TWSCVS}
T.~Sanders, \emph{The $\ell^1$-norm of the Fourier transform on
compact vector spaces}. Preprint.
\bibitem{TWSAS}
T.~Sanders, \emph{Additive structures in sumsets}. Preprint.
\bibitem{TWSZP}
T.~Sanders, \emph{The Littlewood-Gowers problem}. Preprint.

\end{thebibliography}
\end{document}